\documentclass[12pt]{article}

\usepackage{amsmath,amsfonts,amsthm,amssymb, commands,graphicx,listings,enumerate,float,multirow}
\usepackage{geometry,algpseudocode,algorithmicx,algorithm,threeparttable, color,booktabs}
\usepackage{graphicx}
\usepackage{subcaption}
\usepackage[sort]{natbib}
\usepackage{hyperref}

\newcommand{\MIQP}{\operatorname{MIQP}}
\newcommand{\MISOCP}{\operatorname{MISOCP}}
\newcommand{\SOCP}{\operatorname{SOCP}}
\newcommand{\dsum}{\displaystyle \sum}
\newcommand{\R}{\mathbb{R}}
\newcommand{\K}{\mathcal{K}}
\newcommand{\argmin}{\operatorname{argmin}}
\newcommand{\interior}{\operatorname{int}}
\title{Pricing in non-convex markets with quadratic deliverability costs}

\author{Xiaolong Kuang
        \thanks{Department of Industrial and Systems Engineering, Lehigh University, Bethlehem, PA, USA.  E-mail: {\tt xik312@lehigh.edu}}
        \and
               Alberto J. Lamadrid
        \thanks{Department of Economics, Lehigh University, Bethlehem, PA, USA. E-mail: {\tt ajlamadrid@lehigh.edu}}
         \and
         Luis F. Zuluaga
        \thanks{Department of Industrial and Systems Engineering, Lehigh University, Bethlehem, PA, USA.  E-mail: {\tt luis.zuluaga@lehigh.edu}}}

\date{\today}

\begin{document}

\maketitle

\begin{abstract}
The problem of obtaining market-clearing prices for markets with non-convexities has been widely studied in the literature. This is particularly the case in electricity markets, where worldwide deregulation leads to markets in which non-convexities arise from the decisions of market operators regarding which generators are committed to provide electricity power.
Here, we extend seminal results in this area to address the problem of obtaining market-clearing prices for markets in which beyond non-convexities, it is relevant to account for convex quadratic market costs. In a general market, such costs arise from quadratic commodity costs or transactions costs.
 In an electricity market, such quadratic costs arise when ramping costs need to be considered due to the presence of renewable energy sources, which continue to increase their participation in electricity markets. To illustrate our results, we compute and analyze the clearing prices of a classical market problem with the addition of ramping costs.
\end{abstract}

\smallskip
\noindent \textbf{Keywords\ } market-clearing prices, quadratic costs, ramping costs, unit commitment, mixed-integer quadratic programming, renewable energy sources.

\section{Introduction}\label{sec.introduction}

The presence of non-convexities  is inherent in markets due to 
economies of scale, start-up and/or shut-down costs, avoidable costs, indivisibilities, and minimum supply requirements. These non-convexities make the problem of finding appropriate prices that result in a market equilibrium challenging.
This issue has been addressed in classical work by \citet{starr1969quasi,  gomory1960integer, ref.scarf, scarf1990mathematical}. Continued work in this area has been recently reviewed by~\citet{ref.liberopoulos2016critical} and~\citet{van2011}.
Here, we focus on considering the effect of not only non-convexities, but also potential convex quadratic costs that affect market prices. As discussed below, these market characteristics are particularly important in deregulated electricity markets with high penetration of renewable energy sources.

In electricity deregulated markets~\citep[see, e.g.,][]{ref.ruiz2012pricing, ref.oneill2005}, a market operator or individual system operator (ISO) receives  information about generation constraints, marginal generation costs, and fixed commitment costs, from generators participating in the market. Based on this information, the ISO decides the generators that should be committed in the market, as well as their appropriate compensation.

Due to the electricity market non-convexities related to start-up/shut-down costs and minimum output requirements (among others), it is difficult for ISOs to obtain the appropriate compensation values for committed generators. One way in which ISOs can obtain an estimate of the correct compensation values is by finding the shadow (dual) prices associated to the {\em linear programming relaxation}~\citep[cf.,][]{conforti2014} of the {\em mixed-integer linear program} (MILP)~\citep[cf.,][]{conforti2014} associated to the
market's {\em unit commitment problem}~(UC)~\citep[cf.,][]{hobbs2006}. The generators' compensation values are then computed by adding {\em uplift  payments} \citep[cf.,][]{932339}. However, these
uplift payments may be
significant enough to modify the suppliers'
incentives~\citep{ref.liberopoulos2016critical}.
To address this issue, a number of alternative pricing schemes have been recently developed
by \citet{ref.garci2006electricity, ref.oneill2005, bjorndal2008equilibrium, hogan2003minimum, ref.ruiz2012pricing, ref.liberopoulos2016critical}, among others
\citep[see,][for a recent review]{ref.liberopoulos2016critical}. However, to the best of our knowledge, none of these approaches directly consider potential ramping costs in the electricity market. 

The penetration of renewable energy sources such as wind, solar, and wave energy into the electricity market has been steadily increasing. 
The U.S. Energy Information Administration estimates that in the U.S.,
the percentage of energy generated from renewable energy sources has
increased from~9.5\% in~2006 to~13.3\% in~2015. As a result, ISO's
commitment decisions require conventional generators to ramp up or down regularly due to the variability in the power that renewable sources can provide throughout a given day
\citep[see, e.g.,][]{sioshansi2009evaluating}. At the same time, ramping results in substantial wear and tear costs \citep[see, e.g.,][]{lew2012impacts}. 
Therefore, considering these ramping costs in the computation of compensation values is becoming necessary to provide adequate settlement prices for the different generators in the electricity market.

In order to take into account ramping costs in the electricity market, we extend the seminal results introduced by \citet{ref.oneill2005} to obtain appropriate 
prices in markets with non-convexities. \citet{ref.oneill2005} show that after solving the MILP associated with the ISO's UC; the desired clearing prices can be obtained from the shadow prices of the demand and capacity constraints of the linear program resulting from fixing the commitment (binary) variables of the ISO's UC to their optimal value.

Here, we extend this approach by considering ramping costs (or in more generality, quadratic convex costs) into the formulation of the unit commitment problem. Then, using {\em convex optimization} techniques~\citep[cf.,][]{boyd2004convex}, 
we obtain the market-clearing prices in the presence of ramping costs from the {\em dual variables} values (shadow prices) associated with the optimal solution of an appropriate convex optimization problem. That is,
one obtains a set  of  prices  that
  satisfy  the  market  
equilibrium conditions in which suppliers would not want 
to  change  their energy  dispatch  at  these  prices.  
With these results in hand, we perform numerical experiments to show the impact of ramping costs on the clearing prices of a market with non-convexities.

To make the article more self-contained, we introduce some of the basic convex optimization tools used in the article and provide appropriate references for further relevant results in this area.

The pricing results are mainly motivated by an ISO looking to obtain clearing prices in an electricity market with renewable energy sources. However, following~\citet{ref.oneill2005}, we use the general formulation of a bidding market in which the auctioneer is buying and/or selling a commodity, and has an objective of maximizing the value to bidders, when potential quadratic commodity costs~(e.g., in a labor market~\citep{bentolila1992}) or quadratic transaction costs (e.g., in a finance market~\citep{olivares2015, brown2011} are present in the market).

It is worth to mention that related advances on obtaining market-clearing prices in markets with non-convexities have been done recently. For example, consider the work of: \citet{zoltowska2016}, who considers demand shifting bids and transmission constraints in the market; \citet{ye2015}, who consider, from the point of view of the customer, non-convexities arising when considering flexible demand options (e.g., the option to forgo demand); \citet{van2011}, who propose a new model to obtain prices in a market with non-convexities by combining the approaches used in both the US and European electricity markets; and \citet{muatore2008}, who proposes an algorithm to find market-clearing prices in markets with non-convexities that provide incentives to both maintain  and increase generating  capacity.

The rest of the article is organized as follows. In Section~\ref{sec.mainresults}, we introduce the market's assignment problem in the presence of convex quadratic costs. Furthermore, we generalize the results of~\citet{ref.oneill2005} to obtain clearing prices for the market. In Section~\ref{sec.numeric}, we illustrate our results by computing and analyzing the clearing prices associated to Scarf's classical market problem~\citep{ref.scarf} when potential ramping costs are taken into account. In Section~\ref{sec.conclusions}, we finish with some concluding remarks.

 \section{Market problem with convex quadratic costs}\label{sec.mainresults}
 
 Consider a bidding market, where the auctioneer is buying and/or selling a commodity, and has an objective of maximizing the value to bidders, when potential transaction convex quadratic costs are present in the market. Following~\citet{ref.oneill2005}, the market's assignment problem can be formulated as the following {\em mixed-integer quadratic program}~(MIQP) \citep[cf.,][]{lazimy1982}:
 \bequation\label{opt.PMIQP} 
    \begin{array}{llllllll}
  p_{\MIQP}^\ast =   & \min  & \dsum_{k=1}^nc_kx_k+d_kz_k + r_k (x_k-x_k^0)^2 \\
                     & \st   &  \dsum_{k=1}^na_kx_k = b_0,  \\
                     &       &  g_kx_k + h_kz_k \geq b_k, & k=1,\dots,n, \\
                     &       &  z_k \in \{0,1\}, & k=1,\dots,n,\\
                                         &       &  x_k \geq 0, & k=1,\dots,n, 
    \end{array}
\eequation
where for any market bidder $k=1,\dots,n$: $x_k \in \R_+$ represents the units of commodity provided by the bidder; $z_k \in \{0,1\}$ indicates whether the bidder is committed to provide units of commodity; $c_k, d_k \in \R$, are respectively the variable and fixed costs associated with the
bidder's activities; $a_k \in \R$ 
 reflects the production or demand characteristics of the bidder in the market-clearing constraint 
 $\sum_{k=1}^na_kx_k = b_0$, 
where $b_0 \in \R$ is the amount of commodity to be auctioned, with $b_0 \neq 0$ in a one-sided auction and $b_0 = 0$ in a two-sided
auction; $g_k, h_k \in \R$
reflect restrictions on the bidder's operations (e.g., production of a particular plant is limited to
the capacity of that plant); $b_k \in \R$ represents the right hand sides of the internal constraints of the bidder. Also, in an extension of the market assignment problem in~\citet[][Sec. 5]{ref.oneill2005}, $r_k \in \R_+$ denotes the quadratic costs associated with the deviation of the commodity provided by bidder $k=1,\dots,n$ from a target or previous commodity production level $x_k^0 \in \R$. The parameters $r_k \in \R_+, x_k^0 \in \R$ can also be used to model quadratic commodity costs, as well as quadratic transaction costs.

Notice that by assuming that $r_k \ge 0$ for $k=1,\dots,n$, one ensures that~\eqref{opt.PMIQP} has a convex quadratic objective. Also, note that if  $r_k = 0$ for $k=1,\dots,n$, problem~\eqref{opt.PMIQP} is equivalent to the PIP problem in~\citet[][Sec. 5]{ref.oneill2005} when bidders are assumed to provide a single commodity. This single-commodity assumption is made here for ease of presentation and to identify the commodity with power in electricity markets. However, all the results presented thereof extend in straightforward fashion to the more general multi-commodity market considered in~\citet{ref.oneill2005}.

In what follows, we assume that problem~\eqref{opt.PMIQP} is feasible; that is, there is an assignment of the bidders that satisfies both the operating constraints of the bidders as well as the market-clearing constraint. Also, without loss of generality, we assume that $a_k \neq 0$, $k=1,\dots,n$ (i.e., bidders that do not contribute to the market-clearing are not considered),
and $r_k > 0$, $k=1,\dots,n$ (i.e., bidders that do not incur ramping costs do not need a corresponding quadratic constraint in formulation~\eqref{eq.p2} below).

By introducing the auxiliary variables $y_k \in \R_+, k=1,\dots,n$, the market assignment problem~\eqref{opt.PMIQP} can be reformulated as:
\bequation
  \label{eq.p2}
    \begin{array}{lllllllll}
       p_{\MIQP}^\ast =   & \min  & \dsum_{k=1}^n c_kx_k+d_kz_k + r_k y_k \\
                       & \st   &  \dsum_{k=1}^na_kx_k = b_0,  \\
                       &       &  g_kx_k + h_kz_k \geq b_k, & k=1,\dots,n, \\
                       &       &  y_k \geq (x_k-x_k^0)^2, & k=1,\dots,n,\\
                       &       &  z_k \in \{0,1\}, & k=1,\dots,n,\\
                                              &       &  x_k, y_k \geq 0, & k=1,\dots,n.
    \end{array}
\eequation
This follows from the fact that for any optimal solution of~\eqref{eq.p2}, the constraints $y_k \geq (x_k-x_k^0)^2$, $k=1,\dots,n$, will be tight. Now let $\K^{n+1} \subseteq \R^{n+1}$ denote the {\em second-order} (or Lorentz) cone~\citep[cf.,][]{ref.socp} in dimension $n+1$; that is,
\bequation
\label{eq.socpdef}
  \K^{n+1} = \{ (w_0,w)\in\Rmbb^{n+1}: w_0 \geq \|(w_1,w_2,\ldots,w_n)\|_2 \},
\eequation
where $\| \cdot \|_2$ represents the Euclidean norm. Note that like $\R^n$ or $\R^n_+$, the cone $\K^{n+1}$ is a {\em closed convex cone}~\citep[cf.,][]{ben2001}. Moreover, $\R^n$, $\R^n_+$, $\K^{n+1}$ are {\em self-dual} cones; that is, $(\R^n)^* = \R^n$, $(\R^n_+)^* = \R^n_+$, $(\K^{n+1})^* = \K^{n+1}$, where for any convex set, the dual cone of $S \subseteq \R^n$ is equal to  $S^* = \{u \in \R^n: u^\intercal w \ge 0 \text { for all } w \in S\}$. Here, for any $u,w \in \R^n$, $u^\intercal w = \sum_{i=1}^nu_iw_i$ denotes the usual inner product in $\R^n$. The self-duality of  $\R^n$, $\R^n_+$, and~$\K^{n+1}$ will be key in deriving an appropriate dual problem associated with the continuous relaxation 
(i.e., when $z_i \in \{0,1\}$ is replaced by $z_i \in [0,1]$ for $k=1,\dots,n$) of problem~\eqref{eq.p2}.

From the definition~\eqref{eq.socpdef}, it is not difficult to see~\citep[cf.,][eq. (7)]{lobo1997} that for any $k=1,\dots,n$,
\[
y_k \geq (x_k-x_k^0)^2 \Leftrightarrow \left (y_k+1, y_k-1, 2(x_k-x^0_k) \right ) \in \K^{3}.
\]
Thus, problem~\eqref{eq.p2} can be reformulated as the following {\em mixed-integer second-order cone program}~(MISOCP) \citep{terlaky2013}:
\bequation 
 \label{eq.p3}
    \begin{array}{lllllll}
       p_{\MISOCP}^\ast =  & \min  & \dsum_{k=1}^nc_kx_k+d_kz_k +  r_k y_k \\
                      & \st   &  \dsum_{k=1}^na_kx_k = b_0,  \\
                      &       &  g_kx_k + h_kz_k \geq b_k, & k=1,\dots,n, \\
                      &       &  \left (y_k+1, y_k-1, 2(x_k-x^0_k) \right ) \in \K^{3}, & k=1,\dots,n,\\
                      &       &  z_k \in \{0,1\}, & k=1,\dots,n,\\
                                            &       &  x_k,y_k, z_k \geq 0, & k=1,\dots,n.

    \end{array}
\eequation
Like problem~\eqref{eq.p3} above, a MISOCP is an optimization problem that besides a linear objective and linear constraints, has second-order cone constraints, as well as integer or binary variables. Problems with second-order cone constraints are widely used in applications in engineering and science. In particular, this type of constraints appear in  
structural design problems~\citep[][Sec. 3.5]{lobo1997}, electrical engineering~\citep{lubin2016, bienstock2014}, healthcare \citep{mak2014}, and
supply chain management \citep{atamturk2012}.
In many instances, this is a result of the need to take into account the uncertainty of problem
parameters and obtain
solutions that are robust; that is, perform well in different scenarios \cite[see, e.g.,][]{mak2014}.
Moreover the solution of MISOCP problems can be obtained using MISOCP solvers like {\tt MOSEK}, {\tt CPLEX}, and {\tt Gurobi}.

Similar to~\citet{ref.oneill2005}, we next use the optimal solution of~\eqref{eq.p3} to obtain the shadow (dual) prices associated with the market-clearing and bidder operational constraints in~\eqref{eq.p3}. Namely, let
\begin{equation}
z^* = \argmin_{z \in \{0,1\}^n}\{\eqref{eq.p3}\};
\end{equation}
that is, $z^* \in \{0,1\}^n$ is the vector of optimal values of the binary variables in~\eqref{eq.p3}. After replacing $z_k = z^*_k$, $k=1,\dots,n$ in~\eqref{eq.p3}, we obtain the following {\em second-order conic program} (SOCP)~\cite[cf.,][]{lobo1997}:

\bequation\label{opt.psocp}
   \begin{array}{lllllll}
       p_{\SOCP}^\ast =  & \min  & \dsum_{k=1}^nc_kx_k+d_kz_k +  r_k y_k \\
                      & \st   &  \dsum_{k=1}^na_kx_k = b_0,&(p_0)  \\
                      &       &  g_kx_k + h_kz_k \geq b_k, & (q_k) & k=1,\dots,n, \\
                      &       &  \left (y_k+1, y_k-1, 2(x_k-x^0_k) \right ) \in \K^{3}, & (\gamma_k,\alpha_k,\beta_k) & k=1,\dots,n,\\
                      &       &  z_k =z_k^*, & (p_k) & k=1,\dots,n,\\
                                            &       &  x_k,y_k \geq 0, & & k=1,\dots,n.
    \end{array}
\eequation
Like problem~\eqref{opt.psocp} above, a SOCP is an optimization problem that besides a linear objective and linear constraints, has second-order cone constraints. The key characteristic of SOCPs that will be exploited here, is the fact that SOCPs are {\em convex optimization} problems. This follows from the fact that the second-order cone constraint $(w_0,w) \in \K^{n+1}$ is a convex constraint. As a result, there is a rich {\em duality} theory (which generalizes linear programming duality) for these problems, as well as polynomial-time solution algorithms~\citep[cf.,][]{ref.socp}. In turn, these algorithms can be used together with {\em branch \& bound} techniques~\citep[cf.,][]{conforti2014} to solve MISOCP problems like~\eqref{eq.p3}.

Notice that in~\eqref{opt.psocp}, we have associated the dual variables: $p_0$ with the market-clearing constraint, $q_k$ with the $k$-th bidder's operation constraint,
$p_k$ with the $k$-th bidder's commitment constraint, and $(\gamma_k, \alpha_k, \beta_k)$ with the $k$-th second-order constraint, for $k=1,\dots,n$. Also, for ease of notation, let $u\in \R^n$ represent the vector of variables $u_k$, $k=1,\dots,n$. With this notation~\citep[cf.,][]{ref.socp}, the {\em dual} SOCP corresponding to the {\em primal} SOCP problem~\eqref{opt.psocp} can be obtained by constructing the {\em lagrangean dual} of~\eqref{opt.psocp}. Namely, let
\bequation
\label{equ.lagrangian}      
 L(x,y,z,p_0,p,q,\gamma,\alpha,\beta) = \sum_{k=1}^n L_k(x_k,y_k,z_k,p_0,p_k,q_k,\gamma_k,\alpha_k,\beta_k),
 \end{equation}
 where
 \begin{equation}
 \label{eq.lagk}
 \begin{array}{lcl}
 L_k(x_k,\cdots,\beta_k)    &     = & c_kx_k+ d_kz_k +  r_k y_k  + p_0(\frac{1}{k}b_0-a_kx_k) - q_k(g_kx_k+h_kz_k-b_k)\\
&& + p_k(z_k^*-z_k) 
- (y_k+1,y_k-1,2(x_k-x_k^0))^\intercal(\gamma_k,\alpha_k,\beta_k)\\[2ex]
          & = & (c_k-a_kp_0-g_kq_k-2\beta_k)x_k + (d_k-h_kq_k-p_k)z_k\\
          && +(r_k-\gamma_k-\alpha_k)y_k+\frac{1}{k}b_0p_0+b_kq_k+z_k^*p_k-\gamma_k+\alpha_k+2\beta_kx_k^0,
\end{array}
\eequation
for $k=1,\dots,n$, where $p_0\in(\R)^*=\Rmbb$, $q\in(\R^n_+)^* = \Rmbb^n_+$, $p \in (\R^n)^* = \R^n$, and $(\gamma_k,\alpha_k,\beta_k)\in (\K^3)^* = \K^3$, for all $k=1,\dots,n$, are the dual variables or {\em lagrangian multipliers} associated to each of the constraints in problem~\eqref{opt.psocp}. 

From~\eqref{equ.lagrangian} and~\eqref{eq.lagk}, it follows that the dual problem of~\eqref{opt.psocp}: 
\bequation\nn
  \max_{\scriptsize \begin{array}{l} p_0 \in \R, p \in \R^n, q \in \R^n_+,\\ (\gamma_k,\alpha_k,\beta_k) \in \K^3, k=1,\dots,n \end{array}}\ \min_{x\geq0,y\geq0,z\geq0} L(x,y,z,p,q,u,v,w,\gamma,\alpha,\beta),
\eequation
is equivalent to
\bequation\label{opt.dsocp}
\begin{array}{lllllll}
    d_{\SOCP}^\ast =  & \max  & b_0p_0+\dsum_{k=1}^n \left (b_kq_k+z_k^*p_k- \gamma_k+ \alpha_k+2\beta_kx_k^0 \right )\\
          & \st   &  c_k-a_kp_0-g_kq_k-2\beta_k\geq0, & k=1,\dots,n,  \\
          &       &  d_k-h_kq_k-p_k \geq 0,& k=1,\dots,n,  \\
          &       &  r_k-\gamma_k-\alpha_k \geq 0, & k=1,\dots,n,\\
           &       &  (\gamma_k , \alpha_k,\beta_k) \in \K^3, & k=1,\dots,n,\\
          &       &  q_k \geq 0, & k=1,\dots,n.\\
         \end{array}
\eequation

In Proposition~\ref{prop.strong} below, we show that strong duality holds between~\eqref{opt.psocp} and~\eqref{opt.dsocp}, and that their optimal objectives are attained. For this purpose, we introduce, for any set $S \subseteq \R^n$ the notion of its {\em interior}; that is, $\interior(S) = \{ w \in S: \text{ for any } u \in \R^n, \text{ there exists } \epsilon > 0, \text{ such that } w + \epsilon u \in S\}$. 

\bproposition
\label{prop.strong}
Assume that~\eqref{opt.PMIQP} is feasible, and $a_k \neq 0$, $c_k \ge 0,r_k > 0$, for all $k=1,\dots,n$. Then \[p^\ast_{\MIQP}=p^\ast_{\MISOCP}=p^\ast_{\SOCP}=d^\ast_{\SOCP}.\]
\eproposition
\bproof
  From the discussion above, it is clear that problems~\eqref{opt.PMIQP},~\eqref{eq.p3}, and~\eqref{opt.psocp} are equivalent. Therefore, $p^\ast_{\MIQP}=p^\ast_{\MISOCP}=p^\ast_{\SOCP}$. It then remains to show both $p^\ast_{\SOCP}$ and $d^\ast_{\SOCP}$ are attained and strong duality holds between~\eqref{opt.psocp} and~\eqref{opt.dsocp}; that is, $p^\ast_{\SOCP}=d^\ast_{\SOCP}$. For this purpose, we first show that both~\eqref{opt.psocp} and~\eqref{opt.dsocp} are {\em strictly feasible}~\citep[cf.,][]{ref.socp}. 
  
Notice that from the feasibility of~\eqref{opt.PMIQP}, the fact that the feasible set of~\eqref{opt.PMIQP} is closed, and that the objective of~\eqref{opt.PMIQP} is bounded below
by $\min\{0, k\min_{k=1,\dots,n}\{d_k\}\}$, it follows from Weierstrass' Theorem that~\eqref{opt.PMIQP}  has an optimal solution.  
Let $x^* \in \R^n_+$, $z^* \in \{0,1\}^n$ be the optimal solution of~\eqref{opt.PMIQP}, 
and consider any vector $y \in \R^n_+$ such that $y_k>(x_k^*-x_k^0)^2$, for $k=1,\dots,n$. It is easy to see that $(x^*,y,z^*) \in \R^{2n}_+ \times \{0,1\}^n$ is feasible to~\eqref{opt.psocp}. Furthermore, we have that
  \bequation \nn
    y_k+1 > \|( y_k-1,2(x_k^*-x_k^0))\|_2,  
    \eequation
 for $k=1,\dots,n$. Thus, $(x^*,y,z^*) \in \R^{2n}_+ \times \{0,1\}^n$ is a {\em strictly feasible} solution for~\eqref{opt.psocp}; that is,  $(x^*,y,z^*) \in \R^{2n}_+ \times \{0,1\}^n$ is feasible for~\eqref{opt.psocp}, and $(y_k+1, y_k-1, 2(x^*_k-x_k^0)) \in \interior(\K^3)$, $k=1,\dots,n$. Thus, problem~\eqref{opt.psocp} is strictly feasible. Now consider the assignment
  \bequation
  \label{eq.sol}
  \begin{array}{ll}
  p_0 = \min_{k=1,\dots,n} \left \{\frac{c_k-g_k}{a_k} \right \}, \\[2ex]
   (p_k,q_k,\gamma_k,\alpha_k,\beta_k)= \left (d_k-h_k-1,1,\frac{1}{2}r_k,\frac{1}{4}r_k,-\frac{1}{4}r_k \right ),
   \end{array}
  \eequation
  for $k=1,\dots,n$. Clearly, \eqref{eq.sol} is feasible for~\eqref{opt.dsocp}, with
  \bequation\nn
             \gamma_k = \frac{1}{2}r_k >  \frac{1}{2\sqrt{2}} r_k = \|(\alpha_k,\beta_k)\|_2, 
  \eequation
 for $k=1,\dots,n$. That is,  \eqref{eq.sol} is a feasible solution for~\eqref{opt.dsocp}, and $(\gamma_k, \alpha_k, \beta_k) \in \interior(\K^3)$, $k=1,\dots,n$. Thus, problem~\eqref{opt.dsocp} is strictly feasible.
 The result then follows from SOCP duality~\citep[see, e.g.,][Thm. 13]{ref.socp}.
\eproof

Now let us consider the individual problems associated to the bidders in the market.
For each bidder $k=1,\dots,n$, let $t_0$ be the unit commodity price and $t_k$ be the price reflecting the commitment action offered to individual $k$ by the auctioneer. Define the following problem,
\bequation\nn
    \baligned
                     & \min & & c_kx_k+d_kz_k+r_k(x_k-x_k^0)^2-t_0(a_kx_k)-t_kz_k \\
                     & \st  & &  g_kx_k + h_kz_k \geq b_k,  \\
                     &      & &  x_k\geq 0,  \\
                     &      & &  z_k\in\{0,1\};
    \ealigned
\eequation
which similar to problem~\eqref{opt.PMIQP} is equivalent to the following MISOCP:
\bequation\label{opt.pmisocpk}
    \baligned
                     p_{\MISOCP_k}^*(t_0,t_k) = & \min & &  c_kx_k+d_kz_k+r_ky_k-t_0(a_kx_k)-t_kz_k\\
                     & \st  & &  g_kx_k + h_kz_k \geq b_k,  \\
                     &       & & (y_k+1, y_k-1,2(x_k-x_k^0)) \in \K^3,\\
                     &      & &  x_k, y_k\geq 0,  \\
                     &      & &  z_k\in\{0,1\}.
    \ealigned
\eequation
Similar to~\citet{ref.oneill2005}, below we define both the market-clearing prices and associated market-clearing contracts between the auctioneer and the bidders.

\bdefinition
A competitive equilibrium for the market is a set of prices $\{t_0^\ast,t_k^\ast\}$ and allocations $\{x_k^\ast,z_k^\ast\}$, such that
\begin{enumerate}[(a)]
  \item At the prices $\{t_0^\ast,t_k^\ast\}$, the allocations $\{x_k^\ast,z_k^\ast\}$ solve~\eqref{opt.pmisocpk} for all $k=1,\dots,n$;
  \item The market clears: $\sum_{k=1}^na_kx_k^\ast = b_0$.
\end{enumerate}
\edefinition
\bdefinition
    Let $T_k$ be a contract between the auctioneer and bidder $k \in \{1,\dots,n\}$ with the following terms:
    \begin{enumerate}[(a)]
      \item Bidder $k$ operates following $z_k=z^\ast_k,x_k=x_k^\ast$.
      \item Bidder $k$ receives an amount from the auctioneer that is equal to the following payment: $p_0^\ast a_k x_k + p_k^\ast z_k$.
    \end{enumerate}
\edefinition

In what follows, we refer to $T=\{T_k \text{ for } k=1,\ldots,n\}$, to describe the market-clearing contracts between the auctioneer and the bidders. Below, we provide the main result of the article; namely, a characterization of the market-clearing prices for a market with non-convexities and convex quadratic costs.

\btheorem\label{thm.main} Assume that~\eqref{opt.PMIQP} is feasible, and $a_k \neq 0$, $c_k \ge 0,r_k > 0$, for all $k=1,\dots,n$.
    Let $\{x_k^\ast,y_k^\ast,z_k^{\ast}\}$ for all $k=1,\dots,n$ be an optimal solution to~\eqref{opt.PMIQP}. Also, let $p_0^\ast$, and $\{p_k^\ast,q_k^\ast,\gamma_k^\ast,\alpha_k^\ast,\beta_k^\ast\}$ for all $k=1,\dots,n$ be an optimal solution to~\eqref{opt.dsocp}. If in~\eqref{opt.pmisocpk} we define $t_0=p_0^\ast$ and $t_k=p_k^\ast$ for all $k=1,\dots,n$, then the prices $\{p_0^\ast,p_k^\ast\}$ and allocations $\{x_k^\ast,z_k^{\ast}\}$ for all $k=1,\dots,n$ represent a competitive equilibrium.
\etheorem
\bproof
Note that $\{x_k^\ast,y_k^\ast,z_k^{\ast}\}$ for all $k=1,\dots,n$ is also an optimal solution of~\eqref{opt.psocp}. From the KKT conditions associated to the optimal solutions of both~\eqref{opt.psocp} and~\eqref{opt.dsocp}, it follows that:
    \begin{align}
        & 0  \leq (c_k - a_kp_0^{\ast}-g_kq_k^\ast-2\beta_k^\ast) \bot x_k^\ast \geq 0, & k=1,\dots,n, \label{eq.comp1}\\[1.1ex]
        & 0  \leq (d_k - h_kq_k^\ast - p_k^\ast) \bot z_k^\ast \geq 0, & k=1,\dots,n, \nn \\[1.1ex]  
        & 0  \leq (r_k-\gamma_k^\ast-\alpha_k^\ast) \bot y_k^\ast \geq 0, & k=1,\dots,n,  \label{eq.comp3}\\[1.1ex]
        & 0 = p_0^\ast(b_0-\dsum_{k=1}^n a_kx_k^\ast), \nn\\[1.1ex]  
        & 0 \leq q_k^\ast \bot (g_kx_k^\ast+h_kz_k^\ast-b_k)\geq 0,& k=1,\dots,n, \label{eq.comp5}\\[1.1ex]
        & 0 = p_k^\ast(z_k'-z_k^\ast),& k=1,\dots,n, \nn\\[1.1ex] 
        & (y_k^*+1,y_k^*-1,2(x_k^*-x_k^0)) \bot (\gamma_k^*,\alpha_k^*,\beta_k^*), & k=1,\dots,n, \label{eq.comp7}
    \end{align}
where the notation $u\bot w$, for $u,w \in \R$ denotes the complementary slackness between $u$ and~$w$. Now consider the following problem under the contract $T$; that is, each individual bidder~$k$ is offered prices $\{p_0^\ast,p_k^\ast\}$, then each participant $k=1,\dots,n$ will solve~\eqref{opt.pmisocpk} with $t_0 = p_0^*$, $t_k = p_k^*$, $k=1,\dots,n$  to minimize their operation cost. That is, each bidder will solve
\bequation\label{opt.pmisocpk2}
    \baligned
                     p_{\MISOCP_k}^*(p_0^*,p_k^*) = & \min & &  c_kx_k+d_kz_k+r_ky_k-p_0^*(a_kx_k)-p_k^*z_k\\
                     & \st  & &  g_kx_k + h_kz_k \geq b_k,  \\
                     &       & & (y_k+1, y_k-1,2(x_k-x_k^0)) \in \K^3,\\
                     &      & &  x_k, y_k\geq 0,  \\
                     &      & &  z_k\in\{0,1\}.
    \ealigned
\eequation
Clearly, $(x_k,y_k,z_k) = (x_k^*,y_k^*,z_k^*)$ is feasible for~\eqref{opt.pmisocpk2}, and the objective value of this solution, denoted $\hat{p}_{\MISOCP_k}(p_0^*, p_k^*)$, is
\bequation\nn
    \begin{split}
      \hat{p}_{\MISOCP_k}(p_0^*, p_k^*)=  & c_kx_k^\ast+d_kz_k^\ast+r_ky_k^\ast-p_0^\ast a_kx_k^\ast-p_k^\ast z_k^\ast.\\
    \end{split}
\eequation
Using the complementarity equations~\eqref{eq.comp3}, \eqref{eq.comp5}, and \eqref{eq.comp7}, it follows that
\bequation\nn
    \begin{split}
      \hat{p}_{\MISOCP_k}(p_0^*, p_k^*) 
       = & c_kx_k^\ast+d_kz_k^\ast+r_ky_k^\ast-p_0^\ast a_kx_k^\ast-p_k^\ast z_k^\ast - q^\ast_k(g_kx^\ast_k+h_kz^\ast_k-b_k) \\
        & - (r_k-\gamma^\ast_k-\alpha^\ast_k)y^\ast_k -  (y_k^*+1,y_k^*-1,2(x_k^*-x_k^0))^\intercal(\gamma_k^*,\alpha_k^*,\beta_k^*),\\
       = & (c_k - a_kp_0^{\ast}-g_kq_k^\ast-2\beta_k^\ast)x^\ast_k +  (d_k  - p_k^\ast - h_kq_k^\ast)z_k^\ast +y^\ast_k(\gamma_k^\ast+\alpha_k^\ast)\\
         & - (y^\ast_k\gamma_k^\ast+\gamma_k^\ast+y^\ast_k\alpha_k^\ast-\alpha_k^\ast -2\beta_k^\ast x_k^0),\\
       = &  q_k^*b_k- \gamma_k^\ast + \alpha_k^\ast + 2\beta_k^\ast x_k^0.
    \end{split}
\eequation
Next we show that $(x_k,y_k,z_k) = (x_k^*,y_k^*,z_k^*)$ is the optimal solution for problem~\eqref{opt.pmisocpk2}. Let $(x_k,y_k,z_k) \in \R^2_+ \times \{0,1\}$ be a feasible solution of~\eqref{opt.pmisocpk2}. It follows that 
$g_kx_k+h_kz_k-b_k\geq0$, $y_k\geq0$ and $(y_k+1, y_k-1,2(x_k-x_k^0)) \in \K^3$. Therefore
\bequation
\label{eq.add}
\begin{split}
   &  q_k^\ast(g_kx_k+h_kz_k-b_k)\geq 0 \\
   &  y_k(r_k-\gamma^\ast_k -\alpha^\ast_k)\geq 0,\\
   &(y_k+1,y_k-1,2(x_k-x_k^0))^\intercal (\gamma_k^*,\alpha_k^*,\beta_k^*) \geq 0,
\end{split}
\eequation
since the feasibility of $p_0^\ast$, and $\{p_k^\ast,q_k^\ast,\gamma_k^\ast,\alpha_k^\ast,\beta_k^\ast\}$ for all $k=1,\dots,n$, for~\eqref{opt.dsocp}
ensures that $q_k^* \ge 0$, $r_k-\gamma^\ast_k -\alpha^\ast_k \ge 0$, and $(\gamma_k^*,\alpha_k^*,\beta_k^*) \in \K^3 = (\K^3)^*$. Also, from the fact that $p_0^\ast$, and $\{p_k^\ast,q_k^\ast,\gamma_k^\ast,\alpha_k^\ast,\beta_k^\ast\}$ for all $k=1,\dots,n$ is feasible for~\eqref{opt.dsocp}, $x_k, z_k \in \R_+$, it follows that:
\bequation
\label{eq.subtract}
\begin{split}
   &  (c_k - a_kp_0^{\ast}-g_kq_k^\ast-2\beta_k^\ast)x_k \ge 0, \\
   & (d_k - h_kq_k^\ast - p_k^\ast)z_k \ge 0. \\
\end{split}
\eequation
Now let $p_{\MISOCP_k}$ denote the objective value of~\eqref{opt.pmisocpk2} associated with the feasible solution $(x_k,y_k,z_k) \in \R^2_+ \times \{0,1\}$; that is
\bequation\nn
    p_{\MISOCP_k} =  c_kx_k+d_kz_k+r_ky_k-p_0^\ast a_kx_k-p_k^\ast z_k.
\end{equation}    
Using~\eqref{eq.add} and then~\eqref{eq.subtract}, we have   
    \bequation\nn
\begin{array}{lcl}
    p_{\MISOCP_k} & \geq & c_kx_k+d_kz_k+r_ky_k-p_0^\ast a_kx_k-p_k^\ast z_k  - q_k^\ast(g_kx_k+h_kz_k-b_k)\\
                             &  & - y_k(r_k-\gamma^\ast_k-\alpha^\ast_k) - (y_k\gamma_k^\ast+\gamma_k^\ast+y_k\alpha_k^\ast-\alpha_k^\ast + 2\beta_k^\ast x_k -2\beta^\ast_k x_k^0 )\\
                            &  = & (c_k - a_kp_0^{\ast}-g_kq_k^\ast-2\beta_k^\ast)x_k +  (d_k - h_kq_k^\ast - p_k^\ast)z_k +q_k^*b_k -\gamma_k^\ast + \alpha_k^\ast +2\beta^\ast_k x_k^0\\
                          &  \geq &  q_k^*b_k- \gamma_k^\ast + \alpha_k^\ast +2\beta^\ast_k x_k^0\\
                          &  = & \hat{p}_{\MISOCP_k}(p_0^*, p_k^*).
\end{array}
\eequation
This shows that $(x_k^\ast,y_k^\ast,z_k^\ast)$ is the optimal solution for~\eqref{opt.pmisocpk2}. Furthermore, the solution $(x_k^\ast,y_k^\ast,z_k^\ast)$ satisfies the market-clearing condition $\sum_{k=1}^na_kx_k^* = b_0$, therefore, $(x_k^\ast,y_k^\ast,z_k^\ast)$ will provide a market-clearing allocation.
\eproof

\section{Scarf's market instance}
\label{sec.numeric}
As an example of a market with non-convexities, consider a problem proposed by~\citet{ref.scarf}. The objective is to minimize the total cost while satisfying the demand in an electricity market. Two types of plants are available to provide the electricity in the market.
The characteristics of each type of plant, including costs and operational constraints are summarized in Table~\ref{tab.scarfproblem}.
\begin{table}[!htb]
  \centering
  {\footnotesize
    \begin{tabular}{lrr}
      \toprule 
                                 & \multicolumn{1}{c}{(type~1 plant)} & \multicolumn{1}{c}{(type~2 plant)}\\
        Characteristics & \multicolumn{1}{c}{Smokestack}  & \multicolumn{1}{c}{High Tech}  \\ \midrule
        Capacity        &           16.00              &            7.00             \\
        Construction cost &         53.00              &           30.00             \\
        Marginal cost &         3.00            &          2.00            \\
        Average cost at capacity &  6.31    &   6.28  \\
        Total cost at capacity  &  101.00       &  44.00 \\
      \bottomrule
    \end{tabular}}
    \caption{Characteristics of Smokestack and High Tech plants \citep{ref.oneill2005}.}
    \label{tab.scarfproblem}
\end{table}

Scarf's market problem can be formulated as the following
mixed-integer linear program 
\bequation\label{opt.scarf}
    \begin{array}{lllll}
                      \min  &  \dsum_{i=1}^5(3x_{1i}+53z_{1i})+\dsum_{j=1}^{10}(2x_{2j}+30z_{2j})\\
                      \st   &  \dsum_{i=1}^5{x_{1i}} + \dsum_{j=1}^{10}{x_{2j}} = D,  \\
                            &  x_{1i}-16z_{1i}\leq 0, & i=1,\dots,5,\\
                            &  x_{2j}-7z_{2j}\leq 0, & j=1,\dots,10, \\
                            &  x_{1i},x_{2j} \geq 0, & i=1,\dots,5, j=1,\dots,10, \\
                             &  z_{1i},z_{2j}\in\{0,1\}, & i=1,\dots,5, j=1,\dots,10,
    \end{array}
\eequation
where $D$ is the total demand. The market-clearing price for this problem has been studied in~\cite{ref.oneill2005}. Table~\ref{tab.scarfresults} summerizes the optimal solution of \eqref{opt.scarf} for different values of the demand. It is clear that as the demand increases, the number of plants of different types used can change dramatically. For example, when the demand is 56, all type~1 plants are closed and eight (8) type~2 plants are open; however, when the demand is 60, two (2) type~1 plants are open while now only four (4) type~2 plants are open. The market-clearing prices in Table~\ref{tab.scarfdualprice} are obtained from the dual (shadow) prices of the linear program obtained from~\eqref{opt.scarf} by fixing its binary variables to their optimal value~\cite[cf.,][]{ref.oneill2005}. Note that for all the instances with different demand, the market-clearing prices remain the same.
\begin{table}[!htb]
  \centering
    {\footnotesize
    \begin{tabular}{crrcrrc}
      \toprule
                        &  \multicolumn{2}{c}{Units of type} &&  \multicolumn{2}{c}{Unit's output}\\
                        \cmidrule{2-3} \cmidrule{5-6}
        Demand &  1 & 2  &&  1  &  2  & Total Cost\\ \midrule
        56 & 0 & 8 && 0 & 56 & 352 \\
        58 & 1 & 6 && 16 & 42 & 365 \\
        60 & 2 & 4 && 32 & 28 & 378 \\
        62 & 3 & 2 && 48 & 14 & 391 \\
        64 & 4 & 0 && 64 & 0 & 404  \\
        66 & 2 & 5 && 31 & 35 & 419 \\
        68 & 3 & 3 && 47 & 21 & 432 \\
        70 & 0 & 10 && 0 & 70 & 440 \\
      \bottomrule
    \end{tabular}}
    \caption{Optimal solution of Scarf's market problem~\eqref{opt.scarf}~\citep{ref.oneill2005}.}
    \label{tab.scarfresults}
\end{table}

\begin{table}[!tbh]
  \centering
  {\footnotesize
    \begin{tabular}{ccc}
      \toprule
        Commodity Price &  Plant 1 Start-up Price & Plant 2 Start-up Price \\ \midrule
        3 & 53 & 23 \\
      \bottomrule
    \end{tabular}}
    \caption{Market-clearing price of Scarf's problem.}
    \label{tab.scarfdualprice}
\end{table}

Now we consider a modified Scarf problem with quadratic ramping costs,
\bequation\label{opt.scarfmod}
    \begin{array}{lllllllll}
                      \min  &  \dsum_{i=1}^5(3x_{1i}+53z_{1i}+r_{1}(x_{1i}-x_{1i}^0)^2)\\ & +\dsum_{j=1}^{10}(2x_{2j}+30z_{2j}+r_{2}(x_{2j}-x_{2j}^0)^2)\\
                      \st   &  \dsum_{i=1}^5{x_{1i}} + \dsum_{j=1}^{10}{x_{2j}} = D,  \\
                            &  x_{1i}-16z_{1i}\leq 0, & i=1,\dots,5,\\
                            &  x_{2j}-7z_{2j}\leq 0, & j=1,\dots,10,  \\
                            &  x_{1i},x_{2j} \geq 0, & i=1,\dots,5, j=1,\dots,10,\\
                             &  z_{1i},z_{2j}\in\{0,1\}, & i=1,\dots,5, j=1,\dots,10,\\
    \end{array}
\eequation
where $x_{1i}^0,x_{2j}^0$, for all $i=1,\dots, 5$, $j=1,\dots,10$, are set to be the optimal unit's generation obtained after solving~\eqref{opt.scarf} with a demand $D = 55$. The quadratic terms in the objective function can be interpreted as the ramping costs of moving from the original output plan set from~\eqref{opt.scarf} with a demand $D = 55$. It is clear that if we set $D=55$ in~\eqref{opt.scarfmod}, then no matter what the values of $r_{1},r_{2}$ are, problem~\eqref{opt.scarfmod} will have the same optimal solution. 

After setting $r_{1}=r_{2}=0.1$, we can see from Table~\ref{tab.scarfmodresults1} that in most cases, the number of type~1 plants and type~2 plants with full capacity remain the same as a result of the ramping costs (in contrast with Table~\ref{tab.scarfresults}). Note that as the demand increases, a type~2 plant with partial capacity is opened for production. Table~\ref{tab.scarfmoddualprice1} summarizes the market-clearing prices obtained using the results in Section~\ref{sec.mainresults}. In order to maintain a competitive equilibrium, the market-clearing prices vary from case to case. In contrast with the start-up prices obtained by~\citet{ref.oneill2005}, note that from the results of Table~\ref{tab.scarfmoddualprice1}, it follows that the start-up price can differ for type~2 plants producing at full capacity, and type~2 plants producing at partial capacity. Specifically, the start-up price of type~2 plants producing at full capacity changes with the demand, while the start-up price of type-2 plants producing at partial capacity 
remains the same for demands between 56 to 68. The start-up price for closed type~2 plants  is the same as the price for type~2 plants with full capacity.
The unit commodity price varies in a small range but it does not remain the same. The start-up price for all type~1 plants remains mostly equal to 53, but this price is different for some demand levels (e.g., compare $D=60$ and $D=62$). \begin{table}[!tbh]
  \centering
 {\footnotesize
    \begin{tabular}{@{\extracolsep{3pt}}crrrrrrrrrr@{}}
      \toprule
      &  \multicolumn{4}{c}{Type~1 Output} & \multicolumn{4}{c}{Type~2 Output}  & \\ \cline{2-5}  \cline{6-9}
             &  \multicolumn{2}{c}{Partial} &  \multicolumn{2}{c}{Full}  &  \multicolumn{2}{c}{Partial} &  \multicolumn{2}{c}{Full } & \multicolumn{2}{c}{Cost} \\
               \cline{2-3} \cline{4-5} \cline{6-7} \cline{8-9}   \cline{10-11}

  Demand            &  \multicolumn{1}{c}{No.} &  \multicolumn{1}{c}{Prod.} &  \multicolumn{1}{c}{No.} &  \multicolumn{1}{c}{Prod.} &  \multicolumn{1}{c}{No.} &  \multicolumn{1}{c}{Prod.} &  \multicolumn{1}{c}{No.} &  \multicolumn{1}{c}{Prod.} &  \multicolumn{1}{c}{Ramp} &  \multicolumn{1}{c}{Total}\\ \midrule
              56   & 3 & 45.00 & 0 & 0  & 1 & 4.00 & 1 & 7 & 1.90 & 377.90 \\
              58   & 3 & 46.50 & 0 & 0  & 1 & 4.50 & 1 & 7 & 2.10 & 383.60 \\
              60   & 0 & 0 & 3 & 48 & 1 & 5.00 & 1 & 7 & 2.50 & 389.50 \\
              62   & 0 & 0 & 3 & 48 & 1 & 7.00 & 1 & 7 & 4.90 & 395.90 \\
              64   & 3 & 47.40 & 0 & 0  & 2 & 9.60 & 1 & 7 & 4.62 & 429.02\\
              66   & 0 & 0 & 3 & 48 & 2 & 11.00 & 1 & 7 & 6.05 & 435.05\\
              68   & 0 & 0 & 3 & 48 & 2 & 13.00 & 1 & 7 & 8.45 & 441.45\\
              70   & 1 & 15.00 & 3 & 48 & 0 & 0 & 1 & 7 & 22.50 & 467.50\\
       \bottomrule
    \end{tabular}
    }
    \caption{Optimal solution of modified Scarf's problem with $r_{1}=0.1,r_{2}=0.1$.}
    \label{tab.scarfmodresults1}
\end{table}

\begin{table}[!htb]
  \centering
   {\footnotesize
    \begin{tabular}{@{\extracolsep{3pt}}crrrrr@{}}
      \toprule
        & &   \multicolumn{2}{c}{Plant 1 Start-up Price} & \multicolumn{2}{c}{Plant 2 Start-up Price} \\ \cline{3-4} \cline{5-6}
          \multicolumn{1}{c}{Demand} &  \multicolumn{1}{c}{Unit Price}      &   Partial  & Full  \& Closed       &  Partial  & Full  \& Closed           \\             \midrule
        56 & 2.80 & 53.00 & 53.00 & 30 & 24.40 \\
        58 & 2.90 & 53.00 & 53.00  & 30 & 23.70 \\
        60 & 3.00 & 53.00 & 53.00 & 30 & 23.00 \\
        62 & 3.40 & 46.60 & 46.60 & 30 & 20.20  \\
        64 & 2.96 & 53.00 & 53.00 & 30 & 23.28  \\
        66 & 3.10 & 51.40 & 51.40 & 30 & 22.30\\
        68 & 3.30 & 48.20 & 48.20 & 30 & 20.90\\
        70 & 6.00 & 53.00 & 5.00 & 2 & 2.00\\
      \bottomrule
    \end{tabular}
    }
    \caption{Market-clearing price of modified Scarf's problem with $r_{1}=0.1,r_{2}=0.1$.}
    \label{tab.scarfmoddualprice1}
\end{table}
Tables~\ref{tab.scarfmodresults2} and Table~\ref{tab.scarfmoddualprice2}  show that after setting $r_{1}=0.1,r_{2}=0.3$ in~\eqref{opt.scarfmod},  
similar conclusions as in the case $r_{1}=r_{2}=0.1$ can be reached. However, 
with a higher ramping cost on type~2 plants, we can see that with demand values between 64 to 68, instead of operating one more type~2 plant to satisfy the demand, the optimal solution suggests to operate an additional type~1 plant instead.
\begin{table}[!htb]
  \centering
  {\footnotesize
    \begin{tabular}{@{\extracolsep{3pt}}crrrrrrrrrr@{}}
      \toprule
       &  \multicolumn{4}{c}{Type~1 Output} & \multicolumn{4}{c}{Type~2 Output}  & \\ \cline{2-5}  \cline{6-9} 
         &  \multicolumn{2}{c}{Partial } &  \multicolumn{2}{c}{Full }  &  \multicolumn{2}{c}{Partial } &  \multicolumn{2}{c}{Full } & \multicolumn{2}{c}{Cost} \\
   \cline{2-3} \cline{4-5} \cline{6-7} \cline{8-9}   \cline{10-11}
  Demand            &  \multicolumn{1}{c}{No.} &  \multicolumn{1}{c}{Prod.} &  \multicolumn{1}{c}{No.} &  \multicolumn{1}{c}{Prod.} &  \multicolumn{1}{c}{No.} &  \multicolumn{1}{c}{Prod.} &  \multicolumn{1}{c}{No.} &  \multicolumn{1}{c}{Prod.} &  \multicolumn{1}{c}{Ramping} &  \multicolumn{1}{c}{Total}\\ \midrule
              56   & 3 & 47.4 & 0 & 0 & 1 & 1.6 & 1 & 7 & 0.78 & 379.18 \\
              58   & 0 & 0 & 3 & 48 & 1 & 3.0 & 1 & 7 & 2.70 & 385.70 \\
              60   & 0 & 0 & 3 & 48 & 1 & 5.0 & 1 & 7 & 7.50 & 394.50 \\
              62   & 0 & 0 & 3 & 48 & 0 & 0 & 2 & 14 & 14.70 & 405.70 \\
              64   & 1 & 9.0 & 3 & 48 & 0 & 0 & 1 & 7 & 4.62 & 435.10\\
              66   & 1 & 11.0 & 3 & 48 & 0 & 0 & 1 & 7 & 12.10 & 445.10\\
              68   & 1 & 13.0 & 3 & 48 & 0 & 0 & 1 & 7 & 16.90 & 455.90\\
              70   & 1 & 15.0 & 3 & 48 & 0 & 0 & 1 & 7 & 22.50 & 467.50\\
       \bottomrule
    \end{tabular}
    }
    \caption{Optimal solution of modified Scarf's problem with $r_{1}=0.1,r_{2}=0.3$.}
    \label{tab.scarfmodresults2}
\end{table}

\begin{table}[!htb]
  \centering
 {\footnotesize
    \begin{tabular}{@{\extracolsep{3pt}}crrrrr@{}}
      \toprule
        &  &  \multicolumn{4}{c}{Start-up Compensation}\\
        \cline{3-6}
        &  &   \multicolumn{2}{c}{Type~1 Output} & \multicolumn{2}{c}{Type~2 Output} \\ \cline{3-4} \cline{5-6}
              \multicolumn{1}{c}{Demand} &  \multicolumn{1}{c}{Unit Price}        &   Partial  & Full  \& Closed       &  Partial  & Full  \& Closed           \\             \midrule
        56 & 2.96 & 53.00 & 53.00 & 30.00 & 23.28 \\
        58 & 3.80 & 40.20 & 40.20  & 30.00 & 17.40 \\
        60 & 5.00 & 21.00 & 21.00 & 30.00 & 9.00 \\
        62 & 6.20 & 1.80 & 1.80 & 30.00 & 0.60  \\
        64 & 4.80 & 53.00 & 24.20 & 10.40 & 10.40  \\
        66 & 5.20 & 53.00 & 17.80 & 7.60 & 7.60\\
        68 & 5.60 & 53.00 & 11.40 & 4.80 & 4.80\\
        70 & 6.00 & 53.00 & 5.00 & 2.00 & 2.00\\
      \bottomrule
    \end{tabular}
    }
    \caption{Market-clearing price of modified Scarf's problem with $r_{1}=0.1,r_{2}=0.3$.}
    \label{tab.scarfmoddualprice2}
\end{table}

\begin{figure}[ht] 
  \begin{subfigure}[b]{0.5\linewidth}
    \centering
    \includegraphics[width=0.95\linewidth]{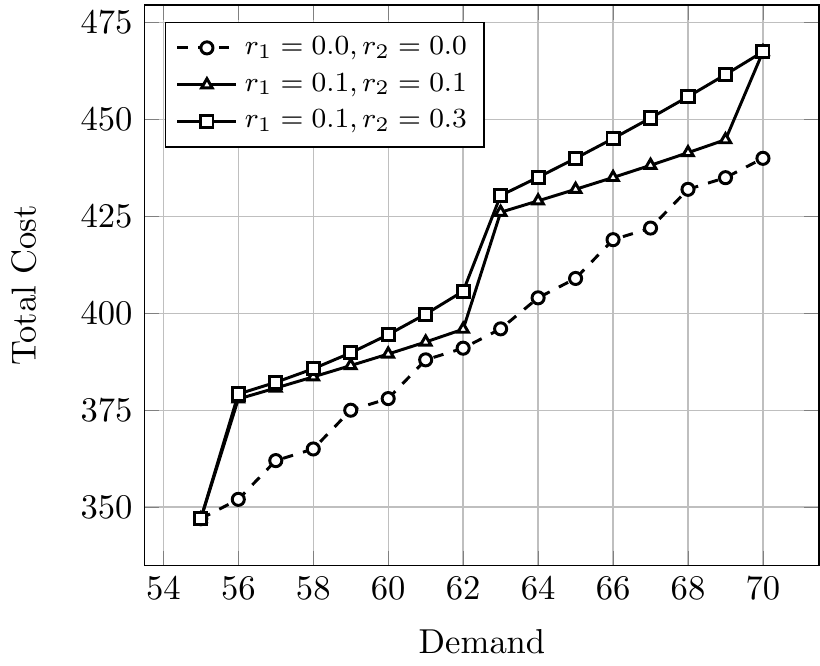} 
    \label{fig7:a} 
    \vspace{4ex}
  \end{subfigure}
  \begin{subfigure}[b]{0.5\linewidth}
    \centering
    \includegraphics[width=0.95\linewidth]{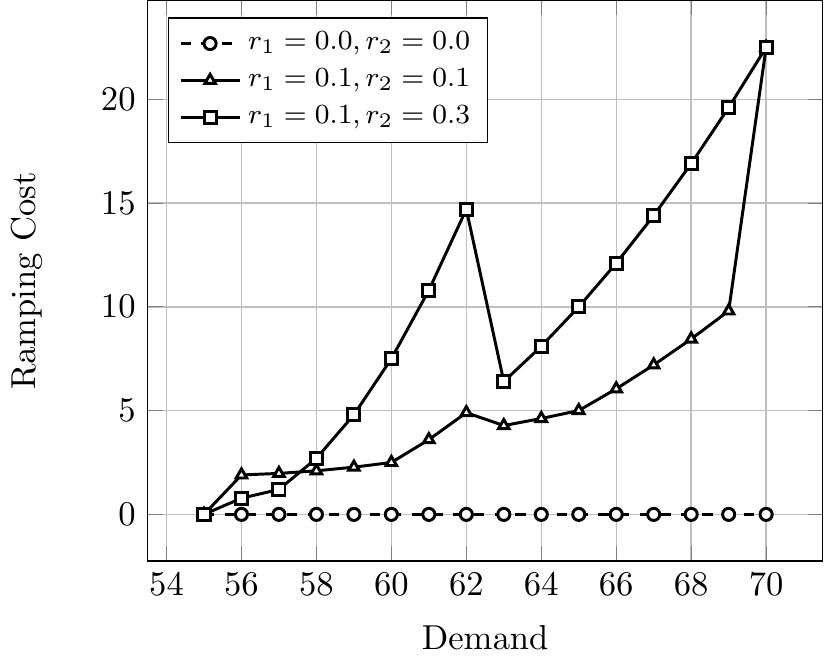} 
    \label{fig7:b} 
    \vspace{4ex}
  \end{subfigure} 
  \begin{subfigure}[b]{0.5\linewidth}
    \centering
    \includegraphics[width=0.95\linewidth]{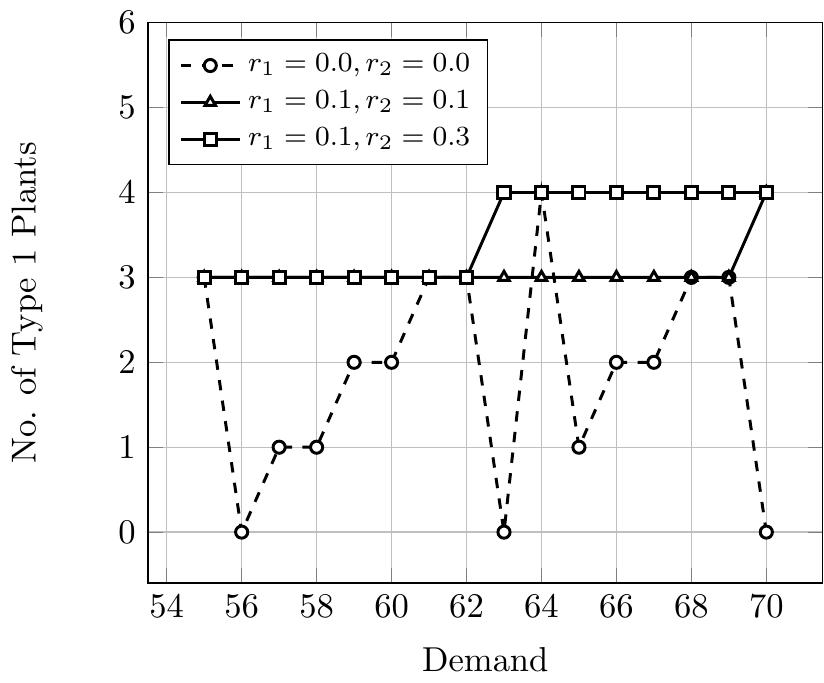} 
    \label{fig7:c} 
  \end{subfigure}
  \begin{subfigure}[b]{0.5\linewidth}
    \centering
    \includegraphics[width=0.95\linewidth]{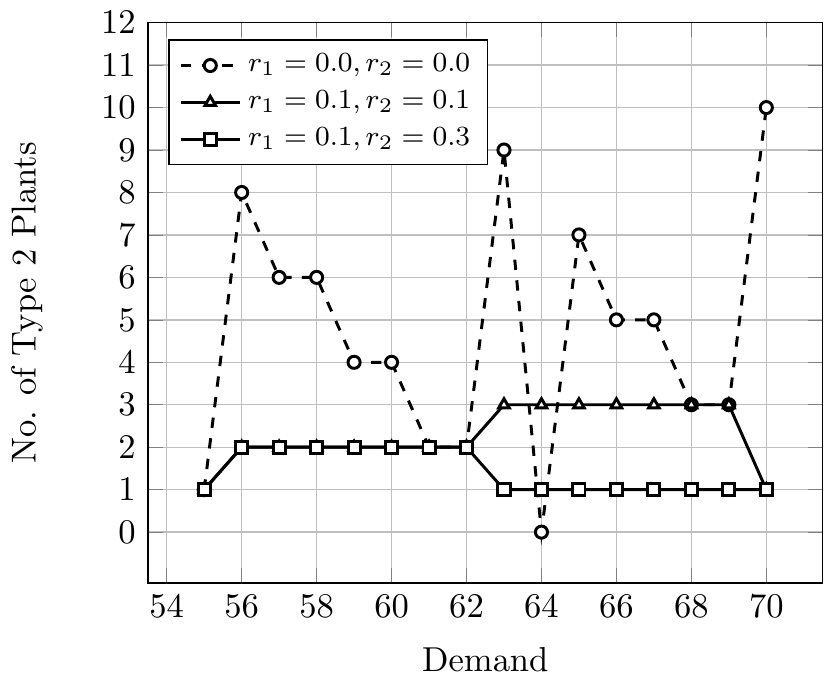} 
    \label{fig7:d} 
  \end{subfigure} 
  \caption{Comparison of optimal solution of~\eqref{opt.scarfmod} for different values of ramping costs $r_1, r_2$ ($r_1 = 0, r_2 = 0$ indicates no ramping costs).}
  \label{fig.optimal} 
\end{figure}

\begin{figure}[!htb]
  \centering
  \includegraphics[width = 0.5\textwidth]{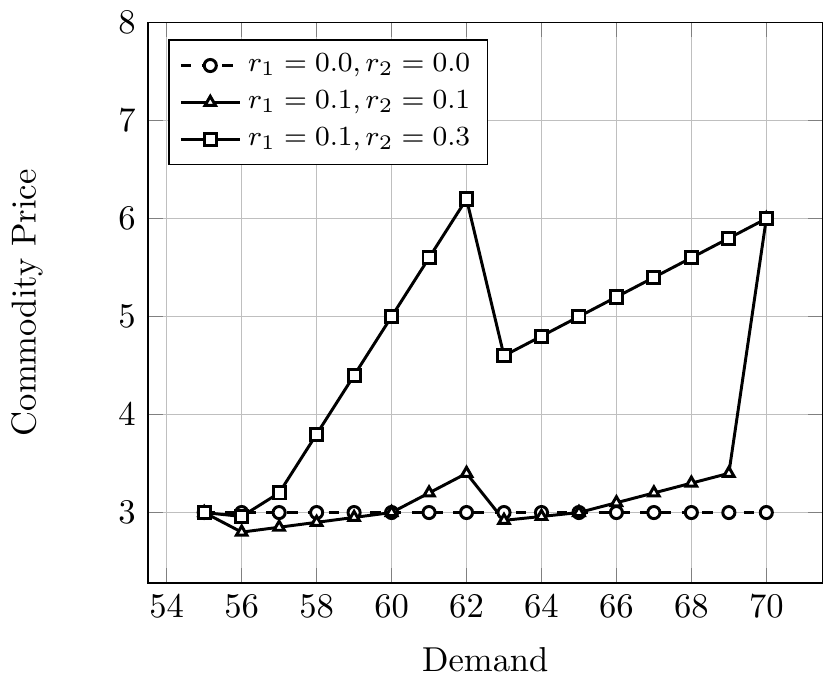}\\
  \caption{Unit commodity prices obtained from the solution of~\eqref{opt.scarfmod}, for different values of ramping costs $r_1, r_2$ ($r_1 = 0, r_2 = 0$ indicates no ramping costs).}\label{fig.unitprice}
\end{figure}

Figure~\ref{fig.optimal} and Figure~\ref{fig.unitprice} graphically compare the solutions obtained from the three cases discussed thus far; that is when $r_1 = r_2 = 0$, when $r_{1}=r_{2}=0.1$, and when $r_{1}=0.1,r_{2}=0.3$.


\begin{table}[!htb]
  \centering
{\footnotesize
    \begin{tabular}{@{\extracolsep{3pt}}crrrrrrrrrr@{}}
      \toprule
      &  \multicolumn{4}{c}{Type~1 Output} & \multicolumn{4}{c}{Type~2 Output}  & \\ \cline{2-5}  \cline{6-9} 
            &  \multicolumn{2}{c}{Partial} &  \multicolumn{2}{c}{Full}  &  \multicolumn{2}{c}{Partial} &  \multicolumn{2}{c}{Full} &  \multicolumn{2}{c}{Cost} \\
      \cline{2-3} \cline{4-5} \cline{6-7} \cline{8-9}   \cline{10-11}
    Demand         &  \multicolumn{1}{c}{No.} &  \multicolumn{1}{c}{Prod.} &  \multicolumn{1}{c}{No.} &  \multicolumn{1}{c}{Prod.} &  \multicolumn{1}{c}{No.} &  \multicolumn{1}{c}{Prod.} &  \multicolumn{1}{c}{No.} &  \multicolumn{1}{c}{Prod.} &  \multicolumn{1}{c}{Ramp} &  \multicolumn{1}{c}{Total}\\ \midrule

        56   & 0 & 0 & 3 & 48 & 1 & 1 & 1 & 7 & 1.0 & 380.0 \\
        58   & 0 & 0 & 3 & 48 & 1 & 3 & 1 & 7 & 9.0 & 392.0 \\
        60   & 0 & 0 & 3 & 48 & 1 & 5 & 1 & 7 & 25.0 & 412.0 \\
        62   & 0 & 0 & 3 & 48 & 1 & 7 & 1 & 7 & 49.0 & 440.0 \\
        64   & 0 & 0 & 3 & 48 & 2 & 9 & 1 & 7 & 40.5 & 465.5  \\
        66   & 0 & 0 & 3 & 48 & 2 & 11 & 1 & 7 & 60.5 & 489.5  \\
        68   & 0 & 0 & 3 & 48 & 2 & 13 & 1 & 7 & 84.5 & 517.5  \\
        70   & 0 & 0 & 3 & 48 & 3 & 15 & 1 & 7 & 75.0 & 542.0  \\
      \bottomrule
    \end{tabular}
    }
    \caption{Optimal solution of modified Scarf's problem with $r_{1}=r_{2}=1$.}
    \label{tab.scarfmodresults3}
\end{table}

Next we consider an extreme case where the parameter of ramping cost is relatively large. With $r_{1}=r_{2}=1$, similar analysis can be applied from the analysis with $r_{1}=r_{2}=0.1$. However, with a relative large ramping cost parameter, the unit commodity price as well as start-up price can change dramatically. It is interesting to see that when the demand becomes 58, the start-up prices for type~1 plant, type~2 plant with full capacity and closed type~2 plant are negative, which means these plants need to pay instead of getting paid to open in exchange for a very high commodity price (compared with Table~\ref{tab.scarfdualprice} and Table~\ref{tab.scarfmoddualprice1}).

Now consider the case in which $r_{1}=r_{2}=1$ and $D=60$ in~\eqref{opt.scarfmod}, as an example to illustrate that the dual prices obtained by using the methodology of Section~\ref{sec.mainresults}, result in a competitive equilibrium. For any $j=1,\dots,n$, the individual problem for a type~2 plant $j$ 
committed by the central operator to produce at
partial capacity is
\bequation\nn
    \baligned
                     & \min & &  2x_{2j}+30z_{2j}+( x_{2j}-0)^2-12 x_{2j}-30z_{2j}\\
                     & \st  & &   x_{2j}-7z_{2j}\leq 0, \\
                     &      & &   x_{2j} \geq 0,  \\
                      &      & &  z_{2j}\in\{0,1\},
    \ealigned
\eequation
and its optimal solution is $(x^\ast_{2j},z^\ast_{2j})=(5,1)$, which matches the optimal solution in Table~\ref{tab.scarfmodresults3}. It is 
not difficult to check, given the dual prices in Table~\ref{tab.scarfmoddualprice3}, that the optimal solution for every individual problem matches its corresponding solution in Table~\ref{tab.scarfmodresults3}, verifying Theorem~\ref{thm.main}.

\begin{table}[!htb]
  \centering
 {\footnotesize
    \begin{tabular}{@{\extracolsep{3pt}}crrrrr@{}}
      \toprule
        & &   \multicolumn{2}{c}{Plant 1 Start-up Price} & \multicolumn{2}{c}{Plant 2 Start-up Price} \\ \cline{3-4} \cline{5-6}
                \multicolumn{1}{c}{Demand} &  \multicolumn{1}{c}{Unit Price}      &   Partial & Full \& Closed       &  Partial & Full \& Closed           \\             \midrule
        56 & 4 & 37 & 37 & 30 & 16 \\
        58 & 8 & -27 & -27 & 30 & -12  \\
        60 & 12 & -91 & -91 & 30 & -40  \\
        62 & 16 & -155 & -155 & 30 & -68 \\
        64 & 11 & -75 & -75 & 30 & -33\\
        66 & 13 & -107 & -107 & 30 & -47\\
        68 & 15 & -139 & -139 & 30 & -61\\
        70 & 12 & -91 & -91 & 30 & -40\\
      \bottomrule
    \end{tabular}
    }
    \caption{Market-clearing price of modified Scarf's problem with $r_{1}=r_{2}=1$.}
    \label{tab.scarfmoddualprice3}
\end{table}

\section{Conclusion}\label{sec.conclusions}
In this paper, we considered the problem of obtaining appropriate market-clearing prices when the market has both non-convexities and convex quadratic costs. Our results show that by using convex optimization techniques, the seminal work of~\cite{ref.oneill2005} on pricing in markets with non-convexity, can be extended to markets in which both non-convexities and convex quadratic costs arise. In the electrical market these two market features arise due to generator fixed costs (or other operational constraints), and ramping or quadratic generation costs. Considering both of these features in the electricity market have become increasingly important due the high penetration of renewable energy sources (RES). This is due to the output volatility of RES, which requires conventional generators to ramp up or down more frequently. Besides electricity markets non-convexity features appear in other markets such as the financial and labor market~\citep{brown2011,bentolila1992}. Thus our results have an impact in a wide range of potential markets.

Furthermore, we believe that the techniques outlined here can be used to extend other pricing methodologies for markets with non-convexities~\citep[see,]{ref.liberopoulos2016critical} that aim to obtain prices with different characteristics. Finally, in the context of electricity markets, it is natural to consider how the consideration of ramping costs will affect prices in an electrical network with congested lines. Addressing these questions provides interesting directions for future research work in this area.

\bibliographystyle{apalike}
\bibliography{references/references}
\end{document}